\documentclass{birkjour}
 \newtheorem{thm}{Theorem}[section]

 \theoremstyle{definition}
 
 \theoremstyle{remark}
 \newtheorem{rem}[thm]{Remark}
 
 \numberwithin{equation}{section}
 \newtheorem{algo}[thm]{Algorithm}

\def\bsymbol#1{%
\mbox{\boldmath$\displaystyle#1$\unboldmath}}

\def\bm#1{\mbox{\boldmath{$#1$}}}
\newcommand{\bfphi}{{\bsymbol{\varphi}}}
\newcommand{\bfv}{{\mathbf{v}}}
\newcommand{\m}{\mathbf{m}}
\def\beq{\begin{equation}}
\def\eeq{\end{equation}}
\def\beqa{\begin{eqnarray*}}
\def\eeqa{\end{eqnarray*}}
\def\beqam{\vspace*{-10mm}\begin{eqnarray}}
\def\eeqam{\end{eqnarray}}
\newcommand{\bfxi}{{\bsymbol{\xi}}}
\newcommand{\bfeta}{{\bsymbol{\eta}}}
\newcommand{\bfzeta}{{\bsymbol{\zeta}}}
\newcommand{\minto}{\textstyle{\int\limits_\Omega}\hspace*{-3.5mm}-}
\newcommand{\mR}{\mathbb{R}}
\newcommand{\pp}{{\mathcal{P}}}
\newcommand{\E}{\mathcal{E}} 
\newcommand{\bfu}{{\mathbf{u}}}

\begin{document}

%
%
%
%
%
%
%
%
%

\title[]
 {Multi-material phase field approach to structural topology optimization}

\author[]{Luise Blank, M.\,Hassan Farshbaf-Shaker, Harald Garcke, Christoph Rupprecht, Vanessa Styles}





\date{\today}

\begin{abstract}

Multi-material structural topology and shape optimization problems are
formulated within a phase field approach. First-order conditions are
stated  and the relation of the necessary conditions to
classical shape derivatives are discussed. An efficient numerical
method based on an $H^1$--gradient projection method is introduced and
finally several numerical results demonstrate the applicability of the
approach. 
\end{abstract}

\maketitle
\section{Introduction}
The efficient use of material and related to that the optimization of shapes and
topology is of high importance for the performance of structures. Many
different methods have been introduced to solve shape and topology
optimization problems and we refer to Bendsoe, Sigmund
\cite{Bendsoe}, Sokolowski, Zolesio \cite{SZ} and Allaire, Jouve, Toader \cite{Allaire2} for
details. In this paper we analyze a multi-phase field approach for
shape and topology optimization problems. This approach is related to
perimeter penalizing methods. However, instead of the perimeter the
Ginzburg-Landau energy 
\begin{align}
E^{\varepsilon}(\boldsymbol{\varphi}):=\int_\Omega \left(\tfrac{\varepsilon}{2} |\nabla{\boldsymbol\varphi}|^2
+\tfrac{1}{\varepsilon}\Psi({\boldsymbol\varphi})\right),\qquad\varepsilon>0,\label{Ginzburg-Landau-Energy}
\end{align}
is added to the objective functional. In
(\ref{Ginzburg-Landau-Energy}) the set $\Omega$ is a given design
domain, the function ${\boldsymbol\varphi}$ which takes values in
$\mathbb{R}^N$ is a phase field vector, $\Psi$ is a potential function
with absolute minima 
which
describe the different materials and
the void and $\varepsilon>0$ is a small parameter related to the
interface thickness. It can be shown that
(\ref{Ginzburg-Landau-Energy}) converges in the sense of
$\Gamma$--limits to the perimeter functional, see
Modica \cite{Modica}. The phase field method
has been introduced in topology optimization by Bourdin and Chambolle
\cite{BC1} and was subsequently used by Burger, Stainko \cite{Burger},
Wang, Zhou \cite{WangZhou}, Takezawa,  Nishiwaki, Kitamura
\cite{Takezawa}, Ded{\'e}, Borden, Hughes \cite{Dedeetal}, Blank et al
\cite{our,BGSSSV} and Penzler, Rumpf, Wirth \cite{Rumpfetal}. However, so
far a rigorous derivation of first order conditions and an analysis of
these conditions in the sharp interface limit $\varepsilon\to 0$ was
missing. In this paper we not only discuss recent progress in this direction
but also introduce and analyze a new efficient method to solve the
constrained minimization problem. 

Although in principle the phase field approach can as well be used for other
shape and topology optimization problems we restrict ourselves
to situations where we seek a domain $\Omega^M$ and a displacement
${\boldsymbol u}$ such that 
\begin{equation}\label{meanc}
\int_{\Omega^M} \boldsymbol{f} \cdot \boldsymbol{u}+\int_{\partial\Omega^{M}} \boldsymbol{g} \cdot \boldsymbol{u}
\end{equation}
or an $L^2$--error to a target displacement
\begin{equation}\label{compl}
\left(\int_{\Omega^M} c |\boldsymbol{u}-\boldsymbol{u}_\Omega|^2\right)^{\frac{1}{2}}
\end{equation}
is minimized subject to the equations of linear elasticity. Here $
\boldsymbol{f}$ and $\boldsymbol{g}$ are volume and surface forces and
$c\ge 0$ is a given weight function on $\Omega$. The optimization problem (\ref{meanc}) is
a mean compliance minimization problem and (\ref{compl}) is an example
of a compliant mechanism problem, see \cite{Allaire2,Bendsoe} for
details. In this contribution we will be brief and refer to \cite{our} and to the 
forthcoming article \cite{BlRu} for details.

\section{Setting of the Problem}\label{Problem}
In this section we introduce how structural topology optimization
problems can be formulated within the phase field approach. 

The goal in multi-material shape and topology optimization is to
partition a given bounded Lipschitz design domain
$\Omega\subset\mathbb{R}^d$ into regions occupied by either void or by
$N-1$ different materials such that a given cost functional is
minimized subject to given constraints. Within the phase field approach
we describe the different material distributions with the help of a
phase field vector $\boldsymbol\varphi:=(\varphi^{i})_{i=1}^{N}$,
where $\varphi^N$ describes the fraction of void and
$\varphi^1, \dots,\varphi^{N-1}$ describe the
fractions of the $N-1$ different materials. The phase field approach
allows for a certain mixing between materials and between materials
and void but the mixing will be restricted to a small interfacial region. In
order to ensure that the phase field vector $\boldsymbol\varphi$
describes fractions we require
that $\boldsymbol \varphi$ lies pointwise in the Gibbs simplex
$\boldsymbol G :=\{ \boldsymbol v \in \mathbb{R}^N \mid  v^i \geq 0\,,\, 
\sum_{i=1}^N v^i =1\}.$ 

%









In this work we prescribe the total
spatial amount of the material fractions through
$\int_{\Omega}\hspace{-5mm}-\hspace{2mm}\boldsymbol\varphi
=\boldsymbol m=(m^i)^N_{i=1}$, where it is assumed that $\sum_{i=1}^N
m^i=1$
with $m^i\in (0,1)$, $i=1,\dots,N$, and where
$\int_{\Omega}\hspace{-5mm}-\hspace{2mm}\boldsymbol\varphi$ denotes
the mean value on $\Omega$. We remark that in principal inequality constraints
for $\int_{\Omega}\hspace{-5mm}-\hspace{2mm}\boldsymbol\varphi$ 
can also be dealt with. 

The potential $\Psi:\mathbb{R}^N\to\mathbb{R}\cup\{\infty\}$ is
assumed to have global minima at the unit vectors $\boldsymbol e_i$,
$i=1,\dots,N$, which correspond  to the different materials 
and  to the void. 

In (\ref{Ginzburg-Landau-Energy}) we  choose an obstacle potential 
$\Psi(\boldsymbol\varphi)=\Psi_{0}(\boldsymbol\varphi)+I_{\boldsymbol
  G}(\boldsymbol\varphi)$ where $\Psi_0$ is smooth and $I_G$ is the
indicator function of the Gibbs-simplex $\boldsymbol G$. Introducing
$\boldsymbol {\mathcal G}:=\{\boldsymbol v\in H^1(\Omega,\mathbb
R^{N})\mid\boldsymbol v(x)\in\boldsymbol G\text{ a.e. in }\Omega\}$
and $\boldsymbol {\mathcal G}^{\boldsymbol m}:=\{\boldsymbol v\in
\boldsymbol {\mathcal
  G}\mid\int_{\Omega}\hspace{-5mm}-\hspace{2mm}\boldsymbol
v=\boldsymbol m\}$ we obtain 
\begin{align}\label{ehat}
\hat{E}^{\varepsilon}(\boldsymbol{\varphi}):=\int_{\Omega}\left(\frac{\varepsilon}{2}|\nabla{\boldsymbol\varphi}|^{2}+\frac{1}{\varepsilon}\Psi_{0}({\boldsymbol\varphi})\right)
\end{align}
and on $\boldsymbol{\mathcal G}$ we have $E^{\varepsilon}(\boldsymbol{\varphi})=\hat{E}^{\varepsilon}(\boldsymbol{\varphi})$. 

We describe the elastic deformation with the help of the displacement
vector $\boldsymbol
u:\Omega\rightarrow\mathbb R^d$ and with the strain tensor
${\mathcal E} ={\mathcal E}(\boldsymbol
u)=
\frac{1}{2}(\nabla\boldsymbol u+(\nabla\boldsymbol u)^T)
$.
The boundary $\partial\Omega$ is divided
into a Dirichlet part $\Gamma_{D}$, a non-homogeneous Neumann part
$\Gamma_{g}$ and a homogeneous Neumann part $\Gamma_{0}$.
Furthermore, 
$\mathbb C$ is the elasticity tensor, $\boldsymbol
f\in L^2(\Omega,\mathbb R^d)$ is the  volume force and
$\boldsymbol g \in
L^2(\Gamma_{g},\mathbb R^d)$ are boundary forces.

The equations of linear elasticity which are  the constraint in our
optimization problem are given by 
\begin{align}\label{state}
\left\{ \begin{array}{rcll}
-\nabla\cdot\left[\mathbb C(\boldsymbol\varphi){\mathcal E}(\boldsymbol u)\right] & = & \left(1-\varphi^{N}\right)\boldsymbol f & \text{ in }\Omega,\\
{\boldsymbol u} & = & {\boldsymbol 0} & \text{ on }\Gamma_{D},\\
\left[\mathbb C(\boldsymbol\varphi){\mathcal E}(\boldsymbol u)\right]{\boldsymbol n} & = & \boldsymbol g&\text{ on }\Gamma_{g},\\
\left[\mathbb C(\boldsymbol\varphi){\mathcal E}(\boldsymbol u)\right]{\boldsymbol n} & = & {\boldsymbol 0}&\text{ on }\Gamma_{0},
\end{array}\right.
\end{align}
where $\boldsymbol n$ is the outer unit normal to
$\partial\Omega$.
The elasticity tensor $\mathbb{C}$ is assumed to depend smoothly on
$\boldsymbol\varphi$, $\mathbb{C}$ has to fulfill the usual symmetry
condition of linear elasticity and has to be
positive definite on symmetric tensors. 
More information and detailed literature on the theory of elasticity
can be found in \cite{our}. 
For the phase field approach the void is approximated by a very soft material with an
elasticity tensor $\mathbb C^N (\varepsilon)$ depending on the interface thickness, e.g.
$\mathbb C^N = \varepsilon^2 \tilde {\mathbb C}^N$ with a fixed tensor 
$ \tilde  {\mathbb C}^N$.
Discussions on how to interpolate the elasticity tensors 
$\mathbb C^i$, for $i=1,\ldots,N$, given in the 
pure materials onto the interface can also be found in Section \ref{landc} 
and in \cite{Bendsoe, our}.

\noindent 
Introducing the notation $\langle{\mathcal A},{\mathcal
  B}\rangle_{\mathbb C}:=\int_{\Omega}{\mathcal A}:\mathbb C{\mathcal
  B}$, where for any matrices ${\mathcal A}$ and ${\mathcal B}$ the
product is given as ${\mathcal A}:{\mathcal
  B}:=\sum_{i,j=1}^{d}{\mathcal A}_{ij}{\mathcal B}_{ij}$, the elastic
boundary value problem (\ref{state}) can be written in the weak
formulation: \\
Given $(\boldsymbol f,\boldsymbol g,\boldsymbol
\varphi)\in L^2(\Omega,\mathbb R^d)\times L^2(\Gamma_{g},\mathbb
R^d)\times L^\infty(\Omega,\mathbb R^N)$ find $\boldsymbol u\in
H^1_D(\Omega,\mathbb R^d)$ such that
\begin{align}
\langle{\mathcal E}(\boldsymbol u),{\mathcal E}(\boldsymbol \eta)\rangle_{\mathbb C(\boldsymbol\varphi)}=\int_{\Omega}\left(1-\varphi^{N}\right)\boldsymbol f\cdot\boldsymbol\eta+\int_{\Gamma_{g}}\boldsymbol g\cdot\boldsymbol\eta=:F(\boldsymbol\eta,\boldsymbol\varphi),\label{mean-compliance}
\end{align}
which has to hold for all $\boldsymbol\eta\in H_{D}^{1}(\Omega,\mathbb R^{d}):=\{\boldsymbol\eta\in H^1(\Omega,\mathbb R^d)\mid\boldsymbol\eta=\boldsymbol 0\text{ on }\Gamma_{D}\}$. The well-posedness of (\ref{mean-compliance}) can be shown by using the Lax-Milgram lemma and Korn's inequality, for details see \cite{our}. \\

\noindent
 Summarized, the structural optimization problem can be formulated as: \\
Given
$(\boldsymbol f,\boldsymbol g,\boldsymbol u_{\Omega}, c)\in
L^2(\Omega,\mathbb R^d)\times L^2(\Gamma_{g},\mathbb R^d)\times
L^2(\Omega,\mathbb R^d)\times L^{\infty}(\Omega)$ and measurable sets
$S_{i}\subseteq\Omega$, $i\in\{0,1\}$, with $S_{0}\cap
S_{1}=\emptyset$, we want to solve
\begin{align*}
(\mathcal P^{\varepsilon})\qquad\begin{cases}
\min & J^{\varepsilon}(\boldsymbol u,\boldsymbol\varphi):=\alpha F(\boldsymbol u,\boldsymbol\varphi)
+  \beta J_{0}(\boldsymbol u,\boldsymbol\varphi)+\gamma \hat{E}^{\varepsilon}(\boldsymbol\varphi),\cr
\text{over } & (\boldsymbol u,\boldsymbol\varphi)\in H_{D}^1(\Omega,\mathbb R^{d})\times H^1(\Omega,\mathbb R^N),\cr
\text{s.t. } & (\ref{mean-compliance}) \text{ is fulfilled and }\boldsymbol\varphi\in \boldsymbol{\mathcal G}^{\boldsymbol m}\cap\boldsymbol U_{c},
\end{cases}
\end{align*}
where $\alpha, \beta \geq 0,\,\gamma, \varepsilon>0$, $\boldsymbol m\in
(0,1)^{N}$ with $ \sum_{i=1}^N m^i=1$, $$\boldsymbol
U_{c}:=\{\boldsymbol\varphi\in H^{1}(\Omega,\mathbb
R^{N})\mid\varphi^{N}=0\text{ a.e. on }S_{0}\text{ and
}\varphi^{N}=1\text{ a.e. on }S_{1}\}$$ and the functional for the
compliant mechanism is given by  
\begin{align}
J_{0}(\boldsymbol u,\boldsymbol\varphi):=\left(\int_{\Omega}\left(1-\varphi^{N}
 \right)\, c\, |\boldsymbol u-\boldsymbol u_{\Omega}|^{2}\right)^{\frac{1}{2}},\label{push-functional}
\end{align}
with a given non-negative weighting factor $c\in L^{\infty}(\Omega)$ 
fulfilling  $|\text{supp }c|>0$. \\
The existence of a minimizer to $(\mathcal P^{\varepsilon})$ 
is shown by classical techniques of the calculus of variations
in \cite{our}. 
\begin{rem}
  From the applicational point of view it might be desirable to fix
  material or void in some regions of the design domain, so the
  condition $\boldsymbol\varphi\in\boldsymbol U_{c}$ makes
  sense. Moreover by choosing $S_{0}$ such that $|S_{0}\cap\text{supp
  }c|\neq 0$ we can ensure that it is not possible to choose only void
  on the support of $c$, i.e. in (\ref{push-functional})
we ensure  $|\text{supp
  }(1-\varphi^{N})\cap\text{supp }c|> 0$.
\end{rem}

\section{Optimality system}\label{OPTIMALITYSYSTEM}
In order to derive first-order necessary optimality conditions for the optimization problem $(\mathcal P^{\varepsilon})$, it is essential to show the differentiability of the control-to-state operator, which is well-defined because of the well-posedness of (\ref{mean-compliance}).  
\begin{thm}\label{Gateaux-diff}
The control-to-state operator $S:L^{\infty}(\Omega,\mathbb R^{N})\rightarrow H_{D}^{1}(\Omega,\mathbb R^{d})$, defined by $S(\boldsymbol\varphi):=\boldsymbol u$, where $\boldsymbol u$ solves (\ref{mean-compliance}), is Fr{\'e}chet differentiable. 
Its directional derivative at $\boldsymbol\varphi\in L^{\infty}(\Omega,\mathbb R^{N})$ in the direction $\boldsymbol h\in L^{\infty}(\Omega,\mathbb R^{N})$ is 
given by $S'(\boldsymbol\varphi)\boldsymbol h=\boldsymbol u^*$, where $\boldsymbol u^*$ denotes the unique solution of the problem
\begin{align}
\langle{\mathcal E}(\boldsymbol u^*),{\mathcal E}(\boldsymbol \eta)\rangle_{\mathbb C(\boldsymbol\varphi)}=-\langle{\mathcal E}(\boldsymbol u),{\mathcal E}(\boldsymbol \eta)\rangle_{\mathbb C'(\boldsymbol\varphi)\boldsymbol h}-\int_{\Omega}h^{N}\boldsymbol f\cdot\boldsymbol\eta,\qquad\forall\boldsymbol\eta\in H_{D}^{1}(\Omega,\mathbb R^{d}).\label{linearized-problem}
\end{align}
\end{thm}
\noindent The expression (\ref{linearized-problem}) formally can be
derived by differentiating the implicit state equation
$\langle{\mathcal E}(S(\boldsymbol\varphi)),{\mathcal E}(\boldsymbol
\eta)\rangle_{\mathbb
  C(\boldsymbol\varphi)}=F(\boldsymbol\eta,\boldsymbol\varphi)$ with
respect to $\boldsymbol\varphi\in L^{\infty}(\Omega,\mathbb
R^{N})$. The proof of Theorem~\ref{Gateaux-diff} can be found in \cite{our}.


\noindent
With Theorem \ref{Gateaux-diff} at hand, we can now derive first order conditions.
Indeed, it follows from the chain rule that the reduced cost functional
$j(\boldsymbol \varphi):=J^{\varepsilon}(S(\boldsymbol
\varphi),\boldsymbol \varphi)$ is Fr\'echet differentiable at every
$\boldsymbol\varphi\in H^1(\Omega,\mathbb R^{N})\cap
L^{\infty}(\Omega,\mathbb R^{N})$ with the Fr\'echet derivative
$j'(\boldsymbol\varphi)\boldsymbol h=J^{\varepsilon}_{'\boldsymbol
  u}(\boldsymbol u,\boldsymbol \varphi)\boldsymbol
u^*+J^{\varepsilon}_{'\boldsymbol\varphi}(\boldsymbol u,\boldsymbol
\varphi)\boldsymbol h$. 
Here we have to assume that $J_0\neq 0$ in case of $\beta\neq 0$. Owing to the convexity of
$\boldsymbol{\mathcal G}^{\boldsymbol m}\cap\boldsymbol U_{c}$, we
have for every minimizer $\boldsymbol\varphi\in\boldsymbol{\mathcal
  G}^{\boldsymbol m}\cap\boldsymbol U_{c}$ of $j$ in
$\boldsymbol{\mathcal G}^{\boldsymbol m}\cap\boldsymbol U_{c}$ that
$j'(\boldsymbol\varphi)(\tilde{\boldsymbol\varphi}-\boldsymbol\varphi)\geq
0, \forall\tilde{\boldsymbol\varphi}\in \boldsymbol{\mathcal
  G}^{\boldsymbol m}\cap\boldsymbol U_{c}$. We can now state the
complete optimality system, see \cite{our} for a proof.
\begin{thm}\label{OPTSYS}
Let $\boldsymbol\varphi\in\boldsymbol{\mathcal G}^{\boldsymbol
  m}\cap\boldsymbol U_{c}$ denote a minimizer of the problem
$(\mathcal P^{\varepsilon})$ and $S(\boldsymbol\varphi)=\boldsymbol
u\in H^1_{D}(\Omega,\mathbb R^d)$, $\boldsymbol p\in
H^1_{D}(\Omega,\mathbb R^d)$  are the corresponding state and adjoint
variables, respectively. Then the functions $(\boldsymbol
u,\boldsymbol\varphi,\boldsymbol p)\in H^{1}_{D}(\Omega,\mathbb
R^d)\times(\boldsymbol{\mathcal G}^{\boldsymbol m}\cap\boldsymbol
U_{c})\times H^{1}_{D}(\Omega,\mathbb R^d)$ fulfill the following
optimality system consisting of the state equation 
\begin{align*}
\text{(SE)}\qquad
\langle{\mathcal E}(\boldsymbol u),{\mathcal E}(\boldsymbol \eta_1)\rangle_{\mathbb C(\boldsymbol\varphi)}=F(\boldsymbol\eta_1,\boldsymbol\varphi),\:
\forall\boldsymbol\eta_1\in H_{D}^{1}(\Omega,\mathbb R^{d}),
\end{align*}
the adjoint equation 
\begin{align*}
\text{(AE)}\:\begin{cases}
\langle{\mathcal E}(\boldsymbol p),{\mathcal E}(\boldsymbol \eta_2)\rangle_{\mathbb C(\boldsymbol\varphi)}\cr
=\alpha F(\boldsymbol\eta_2,\boldsymbol\varphi)+\beta J_{0}^{-1}(\boldsymbol u,\boldsymbol\varphi)\int_\Omega c(1-\varphi^{N})(\boldsymbol u-\boldsymbol u_{\Omega})\cdot\boldsymbol\eta_2 ,\cr
\forall\boldsymbol\eta_2\in H_{D}^{1}(\Omega,\mathbb R^{d}),
\end{cases}
\end{align*}
and the gradient inequality 
\begin{align*}
&\text{(GI)}\:\begin{cases}
\gamma\varepsilon\int_{\Omega}\nabla\boldsymbol\varphi:\nabla(\tilde{\boldsymbol\varphi}-\boldsymbol\varphi)+\frac{\gamma}{\varepsilon}\int_{\Omega}\Psi_{0}'(\boldsymbol\varphi)
\cdot(\tilde{\boldsymbol\varphi}
-\boldsymbol\varphi)\cr
-\frac{\beta}{2} J_{0}^{-1}(\boldsymbol u,\boldsymbol\varphi)\int_{\Omega}c(\tilde{\varphi}^{N}-\varphi^{N})|\boldsymbol u-\boldsymbol u_{\Omega}|^{2}\cr
-\int_{\Omega}(\tilde{\varphi}^{N}-\varphi^{N})\boldsymbol f\cdot(\alpha\boldsymbol u+\boldsymbol p)-\langle{\mathcal E}(\boldsymbol p),{\mathcal E}(\boldsymbol u)\rangle_{\mathbb C'(\boldsymbol\varphi)(\tilde{\boldsymbol\varphi}-\boldsymbol\varphi)}\geq 0,\cr
\forall\tilde{\boldsymbol\varphi}\in \boldsymbol{\mathcal G}^{\boldsymbol m}\cap\boldsymbol U_{c}.
\end{cases}
\end{align*}
\end{thm}

\section{Sharp interface asymptotics}\label{Sharp}
In this section we present the sharp interface limit of the optimality
system given in Theorem \ref{OPTSYS}; for a detailed derivation of the
sharp interface limit using the method of formally matched asymptotic
expansions we refer to \cite{our}. We now consider a more concrete
form of the $\boldsymbol \varphi$-dependent elasticity tensor.
We choose the  
 elasticity tensor starting with
constant elasticity tensors $\mathbb C^{i},\,i\in\{1,\ldots,N-1\}$
which are defined in the pure materials, i.e. when
$\boldsymbol\varphi=\boldsymbol e_{i}$, and model the void as a very
soft material. 
As mentioned, a possible choice of the elasticity tensor in the void is
 $\mathbb C^{N}=\mathbb
C^{N}(\varepsilon)=\varepsilon^2\tilde{\mathbb C}^{N}$ where
$\tilde{\mathbb C}^{N}$ is a fixed elasticity tensor.  In order to
model the elastic properties also in the interfacial region the
elasticity tensor is assumed to be a tensor valued function $\mathbb
C(\boldsymbol\varphi):=\left(\mathbb
  C_{ijkl}(\boldsymbol\varphi)\right)_{i,j,k,l=1}^{d}$ which
interpolate between $\mathbb{C}^1,\dots,\mathbb{C}^{N-1},\mathbb{C}^N(\varepsilon)$.
Furthermore we assume that the weighting factor $c$ in the compliant mechanism functional $J_{0}$ is a smooth function. \\
\noindent
The asymptotic analysis gives that the
phase field functions converge as $\varepsilon $ tends to zero to
a limit function $\boldsymbol \varphi$ which only takes values in
$\{\boldsymbol e_1,...,\boldsymbol e_N\}$. This implies that the
 domain $\Omega$ is partitioned
into $N$ regions $\Omega^{i},\:i\in\{1,\ldots,N\}$, which are
separated by interfaces $\Gamma_{ij},\:i<j$. We choose a unit normal at
$\Gamma_{ij}$ such that for
$\delta>0$ small we have $x+\delta{\boldsymbol\nu}\in \Omega^{j}$ and
$x-\delta{\boldsymbol\nu}\in \Omega^{i}$. Moreover we define
$[\boldsymbol w]_{i}^{j}:=\underset{\delta\searrow
  0}{\lim}\,(\boldsymbol w(x+\delta{\boldsymbol\nu})-\boldsymbol
w(x-\delta{\boldsymbol\nu}))$. We obtain in regions occupied by
material, i.e. for
$i=1,\ldots,N-1$, that the state and the adjoint equation,
respectively, have to hold
\begin{align*}
  -\nabla\cdot\left[\mathbb C^{i}{\mathcal E}(\boldsymbol u)\right] =
  {\boldsymbol f}\text{ and }-\nabla\cdot\left[\mathbb C^{i}{\mathcal
      E}(\boldsymbol p)\right]=\alpha{\boldsymbol f}+\beta
  J_{0}^{-1}(\boldsymbol u,\boldsymbol\varphi)(\boldsymbol
  u-\boldsymbol u_{\Omega})c \,.
\end{align*}
In case of material-material interfaces,
i.e. $\Gamma_{ij}$, $i,j\in\{1,\dots,N-1\}$ we have continuity in the variables $\boldsymbol u$,
$\boldsymbol p$ and continuity for the normal stresses $\mathbb C {\mathcal
  E}(\boldsymbol u)\boldsymbol\nu$ and $\mathbb C {\mathcal
  E}(\boldsymbol p)\boldsymbol\nu$, i.e. for $i,j\in\{1,\ldots,N-1\}$
and $\boldsymbol w\in\{\boldsymbol u,\boldsymbol p\}$ we have
$\left[\boldsymbol w\right]_{i}^{j} = {\boldsymbol 0},\:\left[\mathbb
  C {\mathcal E}(\boldsymbol w)\boldsymbol\nu\right]_{i}^{j} =
{\boldsymbol 0}\: \text{ on }\Gamma_{ij}$. On $\Gamma_{iN}$ we get $\mathbb
C^{i}{\mathcal E}_{i}(\boldsymbol u)\boldsymbol\nu=\mathbb
C^{i}{\mathcal E}_{i}(\boldsymbol p)\boldsymbol\nu=\boldsymbol 0$. 
Moreover we obtain for all $i,j\neq N$
\begin{align}
  0=\gamma \sigma_{ij}\kappa & -\left[\mathbb C {\mathcal
      E}(\boldsymbol u):{\mathcal E}(\boldsymbol p)\right]_{i}^{j}
  +\left[\mathbb C {\mathcal E}(\boldsymbol u)\boldsymbol\nu\cdot (\nabla\boldsymbol p)\boldsymbol\nu\right]_{i}^{j}\nonumber\\
  &+\left[\mathbb C {\mathcal E}(\boldsymbol p)\boldsymbol\nu\cdot
    (\nabla\boldsymbol u)\boldsymbol\nu\right]_{i}^{j}
+\lambda^i-\lambda^j\;\text{ on
  }\Gamma_{ij}\label{Esh}
\end{align}
where $\kappa$ is the mean curvature of $\Gamma_{ij}$ and 
$\boldsymbol\lambda\in \mathbb{R}^N$ are Lagrange multipliers.
We remark that the terms involving $\boldsymbol u$ and $\boldsymbol
p$ generalize the Eshelby traction known from materials science, see
\cite{our}. In addition for all $i\neq N$ it holds
\begin{align*}
  0=\gamma & \sigma_{iN}\kappa+\mathbb C^i{\mathcal E}_{i}(\boldsymbol u) :{\mathcal E}_{i}(\boldsymbol p)-\frac{\beta}{2} J_{0}^{-1}(\boldsymbol u,\boldsymbol\varphi)c\,|\boldsymbol u-\boldsymbol u_{\Omega}|^{2}\\
  &-\boldsymbol f\cdot(\alpha\boldsymbol u+\boldsymbol
  p)+\lambda^{i}-\lambda^{N}\:\text{
    on }\Gamma_{iN}.
\end{align*}
Above the Lagrange multipliers 
$\lambda^{1},\dots,\lambda^N$ sum up to zero and
they
are related to volume constraints 
$\int_{\Omega^i}\hspace{-6mm}-\hspace{2mm} 1=m^i$
which are obtained from the integral
constraints
$\int_{\Omega}\hspace{-4.5mm}-\hspace{2mm}\boldsymbol\varphi=\boldsymbol
m$ in the sharp interface limit.\\
\noindent
In case that void and two or more materials appear junction points
emerge, where e.g. void and two materials meet, see
e.g. Figure~\ref{fig:Fig1}, and it might be desirable in applications
to influence the angles at the junctions. By an appropriate choice of the potential
$\Psi$ the angles at the junctions can be prescribed, see \cite{our}
for details. 

\begin{rem}
  In the case of one material we recover the classical first order
  conditions for the sharp interface structural optimization problem,
  see e.g. Allaire, Jouve, Toader  \cite{Allaire2}. The conditions we derived above generalize
  the first order conditions in \cite{Allaire2} to the multi-phase
  case.
\end{rem}

\section{Numerical methods}
\label{landc}
\subsection{Choice of the potential}
In the previous section we studied  the Ginzburg-Landau energy
with an obstacle potential which leads to an optimization problem with 
inequality constraints. Using instead a smooth potential would lead
to equality constraints  only which are usually easier to handle. 
However, there is a subtle problem, namely, we can not prescribe
 the total spatial amount of the material by
$\int_{\Omega}\hspace{-5mm}-\hspace{2mm}\boldsymbol\varphi
=\boldsymbol m $
since the identification of pure $i$-th-material with $\varphi^i=1$ does not hold 
any longer but the value attained in phase $i$ depends on $\varepsilon$. 
Only in the limit for $\varepsilon \rightarrow 0$ there is a pure $i$-th phase
at $x\in \Omega$ if $\varphi^i(x) =1$.
In Table \ref{smoothy}  the shift of one phase is presented for 
a numerical experiment. The listed values are
the values in areas where the values stay nearly constant,
reflecting a pure phase.
Therefore, the $i$-th material does not have approximately volume $m_i$ by prescribing 
$\minto \varphi^i = m_i$.
Consequently one has to use the obstacle potential or the spatial amount has to be
modelled in a different way.
\begin{table}[h]\centering
\begin{tabular}{|l||r|r|r|r|r|}
 \hline
$\varepsilon$ &  0.02 &  0.01 &   0.005 &  0.0025 &   0.001 \\
\hline
$\varphi^1$  & $\approx$ 1.33942  & $\approx$ 1.21378 & $\approx$ 1.13630 & $\approx$ 1.11450 & $\approx$ 1.05818 \\
  \hline
\end{tabular}
\medskip

\caption{Values for the phase identification 
using the double well potential.}\label{smoothy}
\end{table}
\vspace*{-7mm}
\subsection{Choice of the stiffness tensor on the interface}
The choice of the stiffness tensor on the interface also has a quite severe 
influence on the solution. 
A rough explanation in the presence of one material is the following: 
The stiffest structure has material everywhere.
The mass constraints prohibit this. However, since it is possible to 
choose $\varphi^i \in (0,1) $  on the interface, it can happen that it is
best to have a large mushy region with a mixture of void and material, 
i.e. a broad interface, 
which leads to a stiffer structure.
Therefore the stiffness tensor on the interface should drop down fast but smoothly
from the higher stiffness to the lower stiffness.
We use an quadratic interpolation of the elasticity tensors 
$\mathbb C^1,\ldots,\mathbb C^N$ and set the directional derivative
in direction from the lower to the higher stiffness at the material with the lower
stiffness to zero.
One possibility for $N$-phases is:
\centerline{$
\mathbb C (\boldsymbol \varphi) 
= \sum_{i,j} \mathbb C^{\max\{i,j \}}\varphi^i \varphi^j 
$
}\\
where the tensors are ordered from high to low 
stiffness.
A similar kind of interpolation is used in the SIMP approach for one material and void \cite{Bendsoe}. The choice of the elasticity tensor on the interface influences
also the speed of the numerical algorithm.

\subsection{Projected $H^1$-gradient method}
In this section we focus on the mean compliance problem, i.e. $\beta = 0$
and we use the reduced problem formulation 
$$
	{ \min_{{\mbox{\scriptsize \boldmath{$\varphi$}}} \in \mathcal G^\m }}  \; \; { j(\bfphi)}:=  J^\varepsilon(S(\bfphi)), \bfphi) 
$$
where 
$\mathcal G^\m = \{ \bfxi \in H^1 \mid \ \minto \bfxi = \m,
\ \xi_i \geq { 0}, \ \sum \xi_i\equiv 1 \mbox{ a.e. in } \Omega 
\}
$
is convex and closed and 
{$j$}$: H^1(\Omega,\mR^N)\cap L^\infty(\Omega,\mR^N) \rightarrow \mR $ 
is {Fr\'echet-differentiable}, where
the directional derivatives are given by:
\begin{align} j'(\bfphi)\bfeta
= \gamma\varepsilon {(\nabla \bfphi,\nabla\bfeta)} + 
\frac{\gamma}{\varepsilon}(\Psi_0^\prime(\bfphi),\bfeta) - 
\alpha (\mathbb C'(\bfphi)(\bfeta)
\E(\bfu),\E(\bfu))\; .\label{jprime}
\end{align}
The first-order condition of a general minimization problem
$\min j(\bfphi)$ s.t. $\bfphi \in U$ where $U$ is convex and closed
can be rewritten as a fixed point equation: For any 
$\lambda > 0$ the solution is given as 
$ \bfphi =P_H(\bfphi - \lambda \nabla_H j(\bfphi ))$ 
where $P_H$ is the projection 
onto the convex feasible set $U$ with respect to the scalar product in $H$, see
\cite{GruSa}.
Based on this projected gradient methods have been developed.
We propose to use the following new variant:
\begin{algo}\label{algo1}
Having a current approximation $\bfphi_k$ and  given a positive $\lambda$\\
perform a line-search along the descent direction
$$  \bm v_k:=  {P_H}(\bfphi_k -  \lambda{\nabla_H j(\bfphi_k )}) 
- \bfphi_k$$
to obtain the step length $\beta_k$. Then set
$\bfphi_{k+1}:= \bfphi_k + {\beta_k} \bm v_k$.\\
Stop the iteration if $\| \bm v_{k} \|_{H} < \textrm{tol}$.
\end{algo}
This is not the more known search along the projected gradient path
$\bfphi_{k+1}:= P_H(\bfphi_k - \beta_k \nabla_H j(\bfphi_k ))$
which requires in each line-search step an (expensive) projection.
We can prove a global convergence result \cite{BlRu} which can be found
for convex functions in \cite{GruSa}.
\begin{thm}\label{5.2}
Let $H$ be a Hilbert space, $U\subset H$ be convex, closed and non-empty and
$j:U\to\mR$ be continuously Fr\'echet differentiable.
Then,
every accumulation point $\bfphi^*$  of $\{ \bfphi_k\} $ generated by 
Algorithm \ref{algo1}
is first order critical
if the Armijo step length rule is used.
\end{thm}
The reduced cost functional 
$j$ is differentiable in $H^1(\Omega,\mR^N)\cap L^\infty(\Omega,\mR^N)$,
which is not a Hilbert-space. 
Nevertheless, 
we choose the Hilbert-space 
$ H=\{\bfxi \in { H^1(\Omega, \mR^N)} \mid \minto \bfxi = \bf 0 \} $ 
with the scalar product 
$(\bfxi, \bfeta)_{H}= (\nabla \bfxi , \nabla \bfeta)$.
The gradient does not exist in $H^1$. 
However, since 
\begin{align} 
\tfrac{1}{2}\|(\bfzeta - \bfphi    + \lambda \nabla_H j(\bfphi ) )\|^2_{H} =
\tfrac{1}{2}\|\bfzeta - \bfphi\|^2_{H} + { \lambda j'(\bfphi )(\bfzeta - \bfphi)}  + c
\end{align}
for some constant $c$,
we do not need the $H$-gradient but only the directional derivatives for the
projection.
Hence, we define and use instead of the projection $P_H$ 
the projection type operator $\pp_H$ where $\pp_H(\bfphi,\lambda )$ is given 
by the solution of
\begin{align} 
&\min 
\tfrac{1}{2}\|\bfzeta - \bfphi\|^2_{H} + { \lambda j'(\bfphi )(\bfzeta - \bfphi)} 
\label{proj}
\\
&\mbox{s. t.} \qquad  \minto \bfzeta=\bm m,\quad  \sum_{i=1}^N \zeta^i \equiv 1, \quad
	\zeta^i \geq 0 \;\forall \, i=1,\ldots,N \; . \nonumber
\end{align}
The existence and uniqueness of a solution $\pp_H(\bfphi,\lambda )$ of (\ref{proj}) 
can be shown in our application, see \cite{BlRu}.
Moreover, under some regularity conditions on $j$ which are fulfilled
for our problem, we can 
show the same global convergence result as in Theorem \ref{5.2}, see \cite{BlRu}.
Numerically we solve the obstacle type problem (\ref{proj})  with a primal-dual active
set approach.

\subsection{Scaling}
\begin{figure}[h]
\begin{center}
\includegraphics[width=5.5cm]{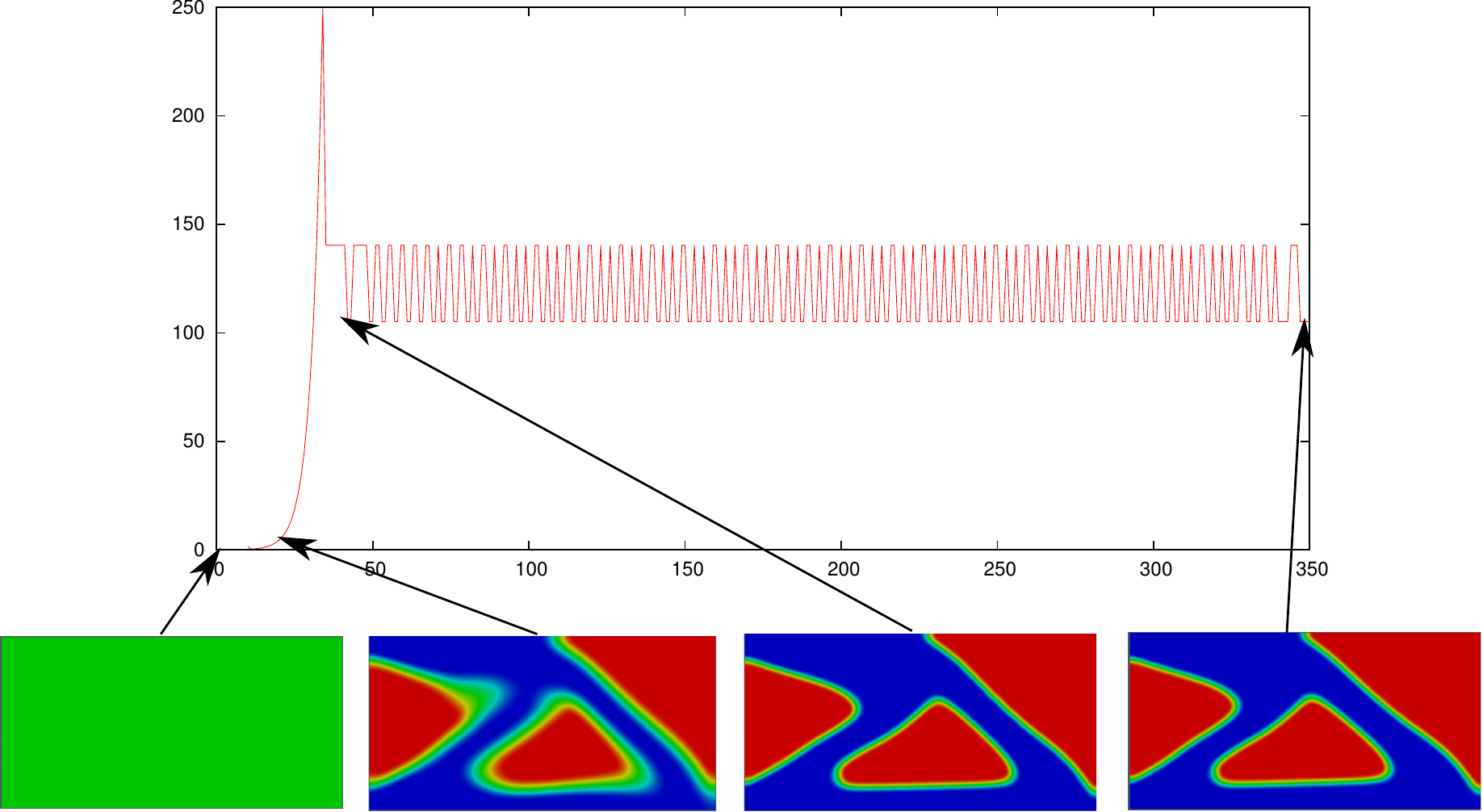}
\caption{Behaviour of $\lambda$ and the phase distribution with respect to the iterations.}\label{lamm}
\end{center}
\end{figure}

In the following we address the choice of the parameter $\lambda $ in the algorithm.
It turned out the scaling of the employed norm is essential for efficiency and for 
iteration numbers independent of the interface thickness, i.e. of $\varepsilon$.
One can motivate this by the fact that the perimeter is approximated 
by the Ginzburg-Landau energy,
which roughly speaking entail 
$\varepsilon || \nabla \bfphi_{\varepsilon}||^2_{L^2} \approx const.$
for the minimizer $\bfphi_{\varepsilon}$.
Hence we have 
$|| \bfphi_{\varepsilon}||_{H} =  O (1/\sqrt{\varepsilon})$
and  $|| \bfphi_{\varepsilon}||_{\sqrt{\varepsilon} H} = O (1)$. 
This is confirmed also 
numerically. As a consequence we choose the ${\sqrt{\varepsilon} H} $ metric. 
Since $\pp_{\sqrt{\varepsilon} H} = \pp_H$ this leads to the use 
of a scaled $H$-gradient since
$\nabla_{\sqrt{\varepsilon} H} \bfphi_{\varepsilon} = 
{\frac1\varepsilon} \nabla_{H} \bfphi_{\varepsilon}$,
respectively this emphasizes to use $\lambda = {\frac1\varepsilon}$.
However, the iterates $\varphi_k$ fulfill $|| \varphi_k||_H \approx || \varphi_\varepsilon ||_H $ only when phases are separated and interfaces are present with thickness according
to $\varepsilon$. In the first iterations this is in general not the case.
Hence, it is more appropriate to adapt $\lambda $ during the iterations.
As a first approach we used the following updating strategy:
\\
{\em
\hspace*{5mm} Set $\lambda_0= \frac{0.01}\varepsilon $ and  choose some  $0<\bar c<1$,
\\
\hspace*{5mm}
in the following  set 
$\lambda_k = \lambda_{k-1} /\bar c $ if $ \alpha_{k-1}=1$  
and $\lambda_k = \lambda_{k-1} \bar c $ else.}
\\
The changes in $\lambda$ with respect to the iterations can be seen examplarily in
Figure \ref{lamm}, where underneath the evolution of the phases can be seen.
We remark that this is no line search with respect to $\lambda$.
In the following algorithm we outline one iteration step and indicate
with it the cost of the method.
\begin{algo}\label{algo2}
\hangindent6mm
Given $\bfphi_k$ and a fixed
$\sigma \in (0,1)$ \\
- solve the { elasticity equation} (\ref{mean-compliance}) for 
$\bfu_k = S(\bfphi_k):\Omega \rightarrow \mR^d$, \\
-  assemble the {directional derivatives} 
$j'(\bfphi_k) \bfeta$\quad $ \forall  \bfeta \in H \cap { L^\infty}  $, \\
- update $\lambda_k$, \\
-  solve the obstacle  type problem (\ref{proj}) for the 
 $\pp_{H^1}(\bfphi_k,\lambda_k ): \Omega \rightarrow \mR^N$,\\
- set $\bm v_k:= \pp_{H^1}(\bfphi_k,\lambda_k ) - \bfphi_k $ and 
stop if $\| \bm v_{k} \|_{\sqrt{\varepsilon} H} < \textrm{tol}$,\\
- determine the {Armijo-step} length $\beta_k= \sigma^{m_k}$ using back tracking\\
\hspace*{5mm}where in each iteration we have to
solve the elasticity equation 
\\ \hspace*{5mm} for 
$\bfu = S(\bfphi_k+ \beta \bfv_k):\Omega \rightarrow \mR^d$, \\
- set $\bfphi_{k+1}:= \bfphi_k + {\beta_k} \bm v_k$ .
\end{algo}

\subsection{Numerical experiments}
\begin{figure}[h]
\begin{center}
\includegraphics[width=2.5cm]{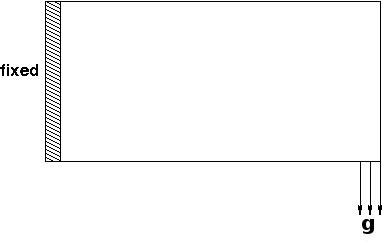}
\hspace*{15mm}
\includegraphics[width=4cm]{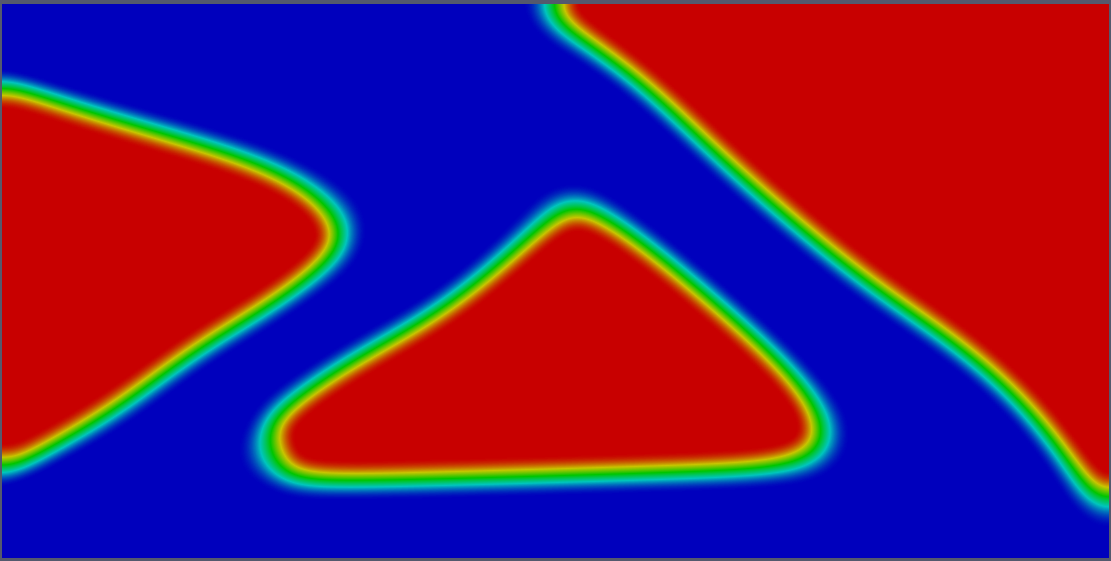}
\caption{Cantileaver beam, geometry (left) and numerical result (right).}\label{clbeam}
\end{center}
\end{figure}

The numerical experiments which underline the above statements are for
the cantilever beam in two dimensions and with one material and void.
The design domain is $\Omega = (-1,1)\times (0,1)$ and $\alpha = 1$. 
There is no volume force  but a boundary force ${\bm g} \equiv (0,-250)^T$ is acting on 
$\Gamma_{g} = (0.75,1)\times \{ 0\}$. The Dirichlet part is 
$\Gamma_D= \{-1\} \times (0,1)$. For the stiffness tensor of the material we take 
$\mathbb C^1\E = 2\mu\E + \lambda (\mbox{tr} \E)I$ with Lam\'e constants $\mu = \lambda =5000$.
Moreover we use the constant $\gamma = 0.5$ and prescribe the masses by $50\%$ material and $50\%$ void. 
Figure \ref{clbeam} displays the
setting and the result for $\varepsilon = 0.03$.

\begin{figure}[h]
\parbox[h]{120mm}{\begin{center}
\hspace*{-10mm}
\begin{minipage}[b]{68mm}%
{\centerline{\includegraphics[width=3.5cm]{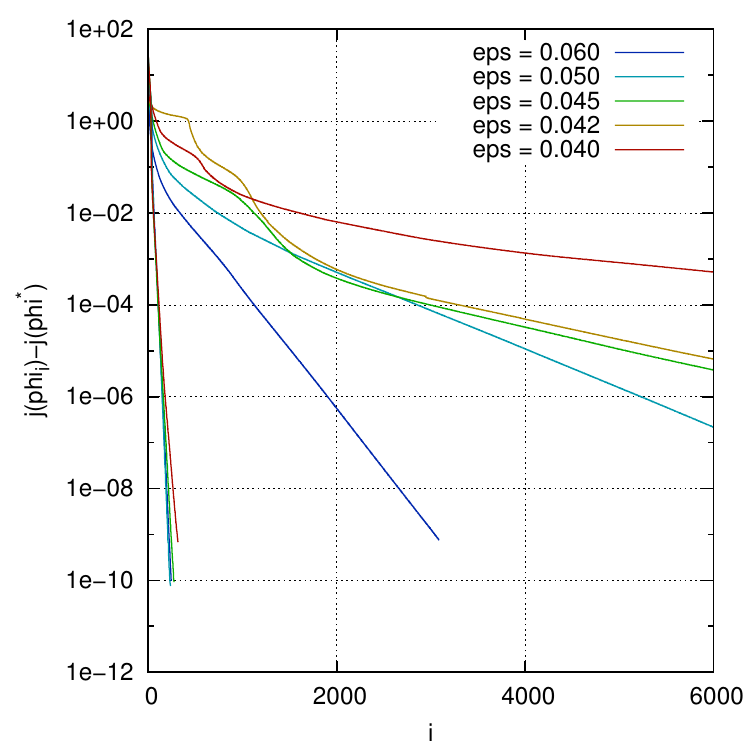}}
\caption{With and without scaling}\label{scaling}}
\end{minipage}
\hspace*{-12mm}
\begin{minipage}[b]{70mm}%
{\centerline{\includegraphics[width=3.5cm]{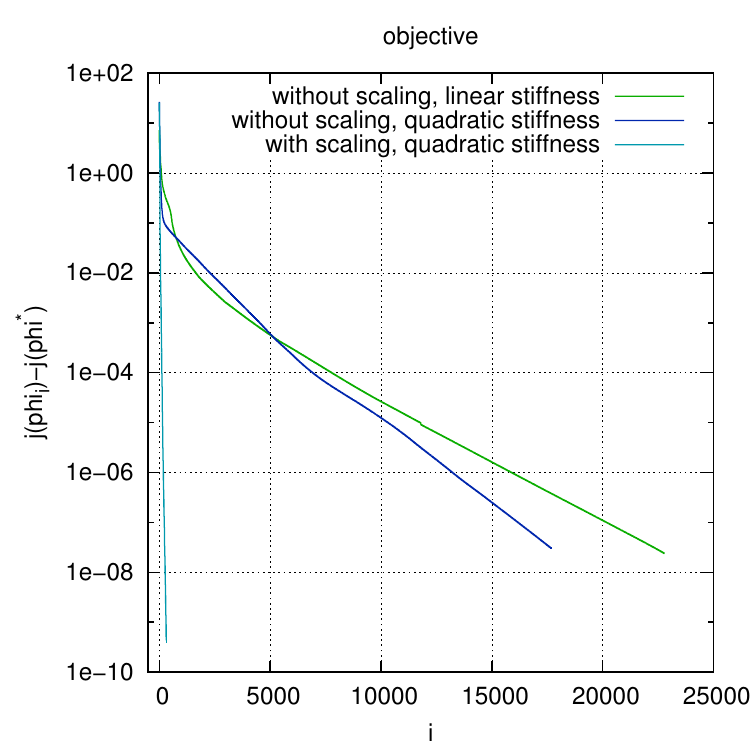}}
\caption{Interpolated elasticity tensor}\label{tens}}
\end{minipage}
\end{center}}
\end{figure}

All computations are done using the finite element toolbox FEniCS \cite{LoggMardalEtAl2012a}.
So far we only use equidistant meshes.
The elasticity equation is discretized with P1-finite elements and the
arising linear systems are solved directly. 
In the computations with one material the problem
setting is reduced to one phase field only by working with 
$\varphi:=\varphi^2-\varphi^1$. In Figure \ref{scaling} the upper
five lines correspond to the results without scaling 
the gradient and shows the approximated error in the cost functional 
with respect to the iteration numbers. 
We clearly see a dependency on 
$\varepsilon$. The lower five lines correspond to the results with
scaling and  are nearly not distinguishable,
independent of $\varepsilon$ and lead to much better approximations
for a lower number of iteration. 
In Figure \ref{tens} the influence
of the choice of the linear versus the quadratic interpolation of the
stiffness tensor is depicted for $\varepsilon= 0.04$.
%

In Table \ref{comp} we study the dependency on the mesh size $h$
and compare the approaches without scaled gradient 
and with linear interpolation of the elasticity tensors (called {\em old}
in the table) to the approach using the scaled gradient and the
quadratic interpolated elasticity tensor (called {\em new}
in the table). In the last column   
we listed the result for the latter approach but using in addition 
nested iteration, i.e. using the result of the previous $h$ as initial
data for the next and solving with an decreasing tolerance $tol$. 
This leads to the expected speed up, here the nested approach needs roughly 
15\% of the CPU-time of the {\em new} approach. 
The more severe speed up of the {\em old} approach is obtained 
by the {\em new} ansatz, which leads to a reduction to 0.5\% 
of the corresponding CPU-time of the {\em old} approach.
Nevertheless, in any case the expected mesh independent number of iterations
is confirmed.
\begin{table}[h]\centering
\begin{tabular}{|l||r|r|r|r|r|r|r|}
  \hline
  & & \multicolumn{2}{c|}{old} & \multicolumn{2}{c|}{new}& \multicolumn{2}{c|}{nested}\\
  h & DOF & CPU & iter. & CPU & iter.  &  CPU & iter.\\
  \hline\hline
$2^{-4}$ & 561   &       12m      &   9956  &    5s &  112&   4s & 85 \\
$2^{-5}$ & 2145  &    2h 25m      &  14590  &    1m &  408&   7s & 52 \\
$2^{-6}$ & 8385  &   20h 40m      &  16936  &    4m &  321&  14s & 24 \\
$2^{-7}$ & 33153 &3d 20h 28m      &  19416  &   21m &  276&   2m & 33 \\
$2^{-8}$ & 131841&{\bf 23d 15h 0m}& 18891&{\bf  3h 1m }&  270&25m &  63\\
         &       &           &         &       &     & total \bf 28m & \\
\hline
\end{tabular}
\medskip

\hspace*{-15mm}
\parbox{13.5cm}{\caption{Comparison of the previous and the new approach as well as with nested iteration for $\varepsilon = 0.04$.}\label{comp}}
\end{table}
We do not list but would like to mention that in the above example
the number of line search iterations stay also mesh independent and are between 1 and 3.
The number of PDAS iterations are mildly mesh dependent 
but stay below 10 after the first few iterations.

\begin{figure}[h]
\begin{center}
\begin{tabular}{c c c }
&Initial data&\\
$\bfphi_0 \equiv${\tiny (38.43,21.33,40.24)}& & random data\\
\includegraphics[width=2.5cm]{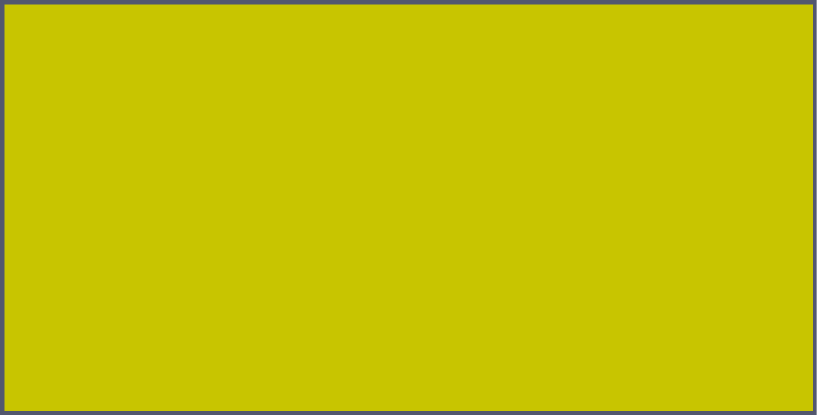}\quad&
\includegraphics[width=2.5cm]{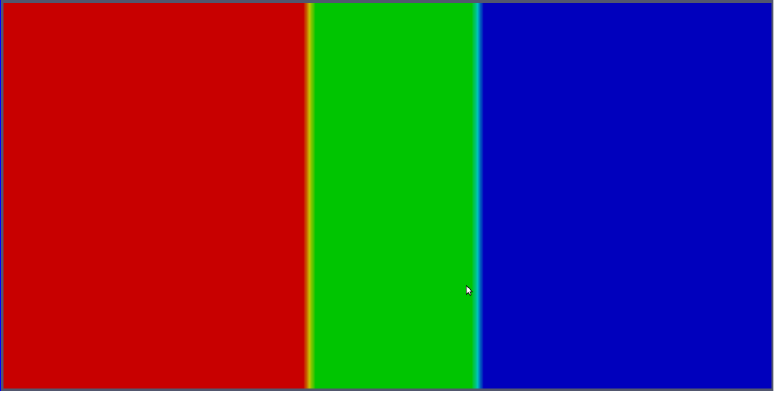}\quad&
\includegraphics[width=2.5cm]{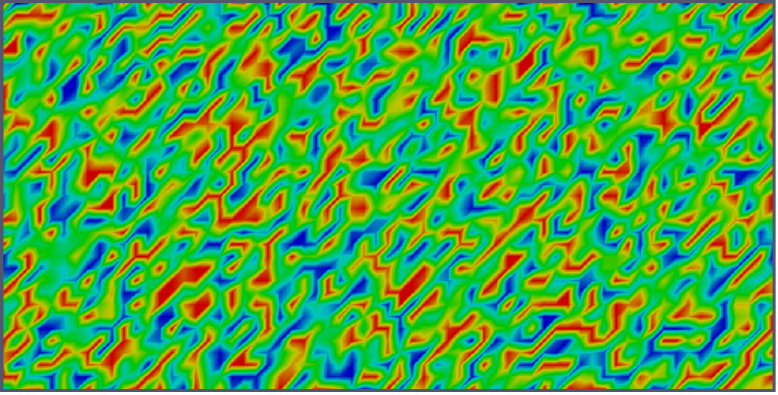}
\\
&Local Minima&\\
\includegraphics[width=2.5cm]{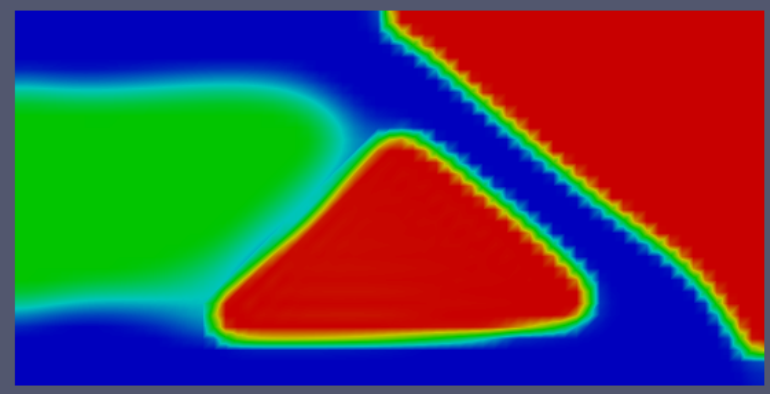}\quad&
\includegraphics[width=2.5cm]{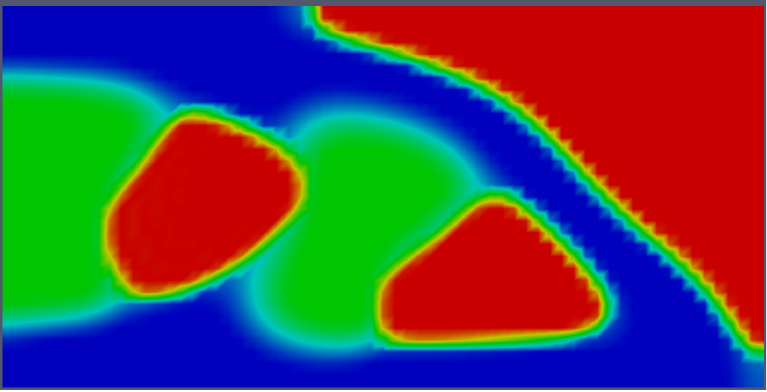}\quad&
\includegraphics[width=2.5cm]{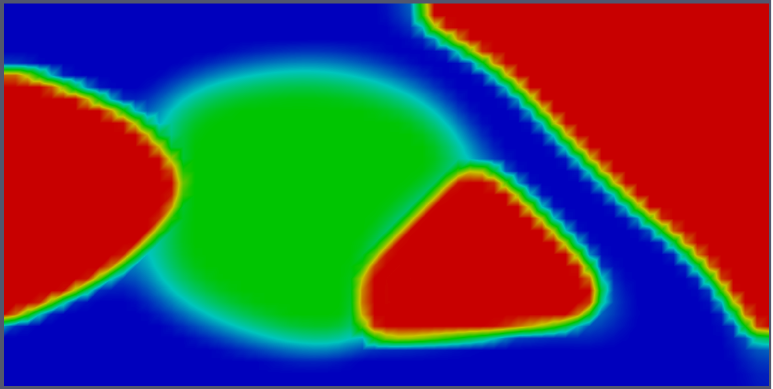} \\
\hspace*{-15mm}$j(\bfphi^*) = $19.4871&19.9138 &18.5631 
\end{tabular}\\
\caption{Cantilever beam with two materials and void in  2$d$.}\label{locmin}
\end{center}
\end{figure}

As expected we can obtain different local minima if we start with different 
initial data as can be seen in Figure \ref{locmin} for a cantilever beam with
two materials and void. The first column shows the result where the initial 
data is a constant mixture of materials and void, the second started with
separated material distribution and the third with random data. The last
yields the lowest value of the cost functional. 

\begin{figure}[h]
\center
\includegraphics[width=6cm]{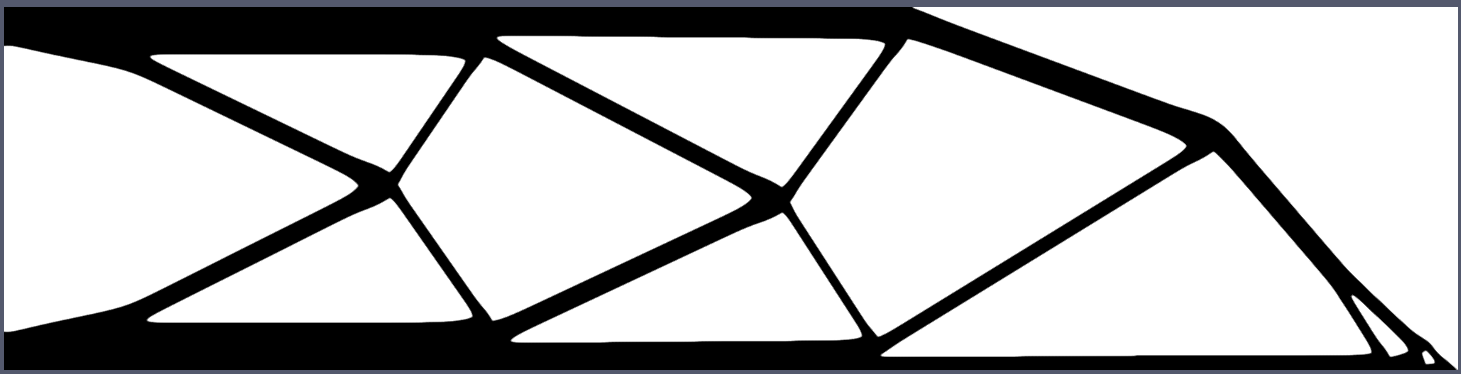}
\caption{A long cantilever beam with low material fraction and a low
  interfacial energy penalization.}
\label{fig:longbeam}
\end{figure}

\begin{figure}[h]
\parbox[h]{120mm}{\begin{center}
\hspace*{-10mm}
\begin{minipage}[b]{68mm}%
{\centerline{\includegraphics[width=3.5cm]{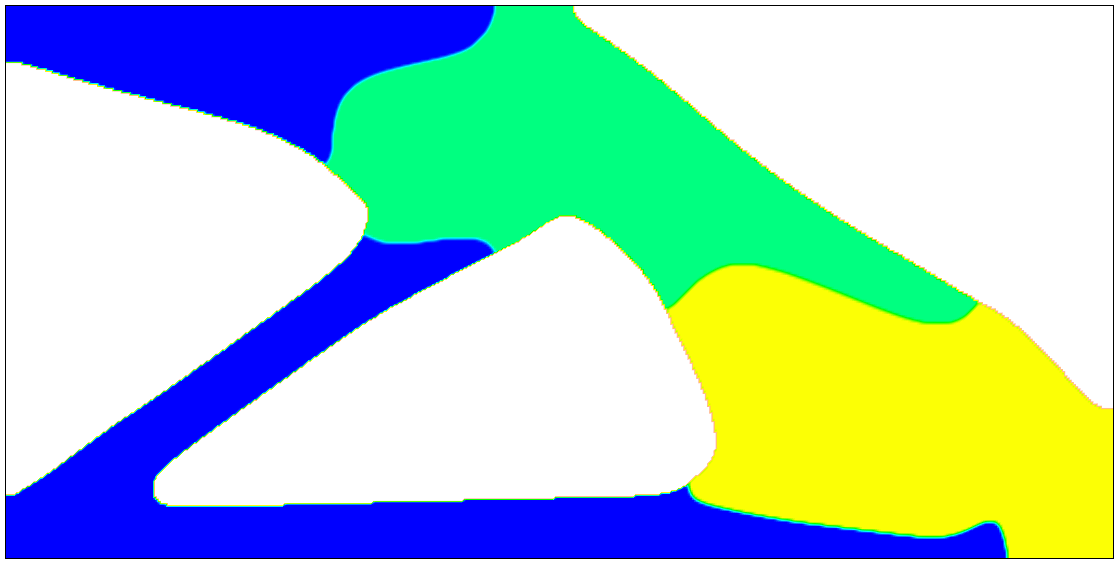}}
\caption{A cantilever beam with four phases.}
\label{fig:Fig1}}
\end{minipage}
\hspace*{-14mm}
\begin{minipage}[b]{72mm}%
{\centerline{\includegraphics[width=3.8cm
]{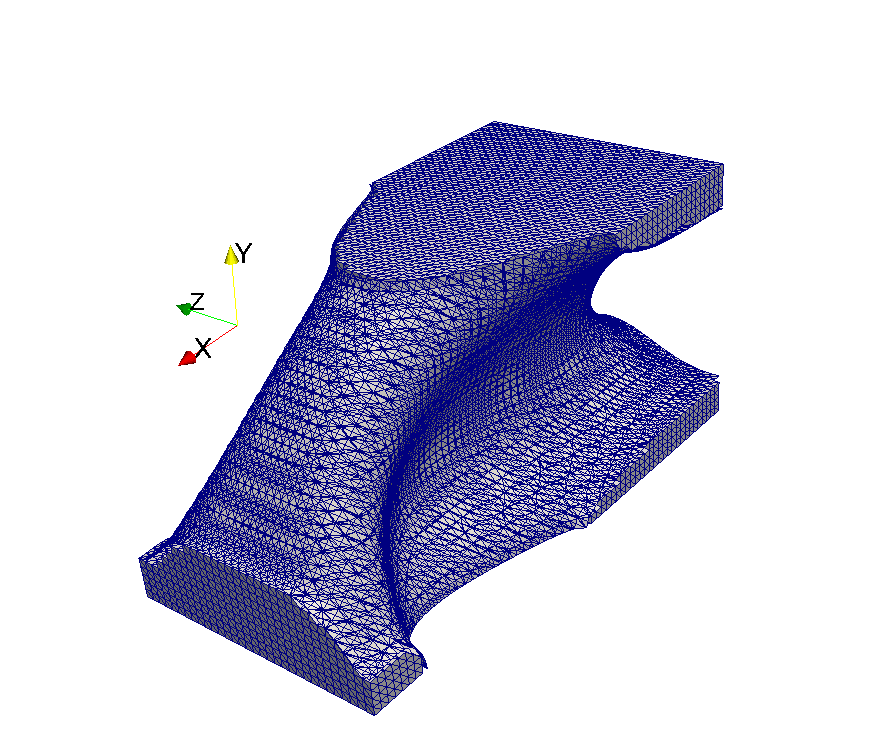}}
\vspace*{-5mm}
\caption{A cantilever beam in three space dimensions.}
\label{fig:Fig2}}
\end{minipage}
\end{center}}
\end{figure}

The following three Figures \ref{fig:longbeam}-\ref{fig:Fig2}  illustrate some results for a long cantilever beam
with one material, for a case with three materials and void and an example
for a cantilever beam in 3d with one material.%

%

\subsection{Numerical results for a compliant mechanism}
\begin{figure}[h]
\centering
\includegraphics[width =3.5cm]{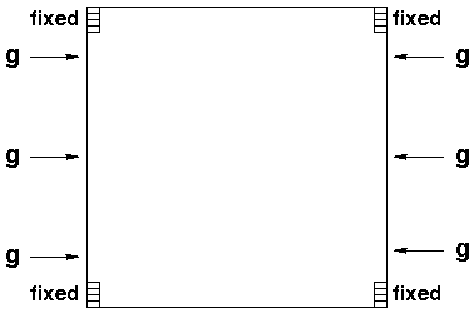}
\vspace*{-5mm}
\caption{Push configuration.}
\label{push_pic}
\end{figure}

In this section we present a compliant mechanism simulation, 
in particular we set $\alpha=0$ in $(\mathcal P^{\varepsilon})$. 
The configuration we consider is depicted in Figure \ref{push_pic}, where zero 
Dirichlet boundary conditions are posed on the left and right boundaries at the top and bottom 
and horizontal forces are applied at sections along the left and right boundaries.

In order to solve the gradient inequality (GI) in Theorem \ref{OPTSYS}, we use 
here as a first numerical approach a classical $L^2$-gradient flow dynamic for 
the reduced cost functional. 
The gradient flow yields the following parabolic variational inequality 
for all $\tilde{\boldsymbol\varphi}\in \boldsymbol{\mathcal G^m}$ and all $t>0$:
\begin{eqnarray}
\varepsilon \int_{\Omega}\frac{\partial \boldsymbol\varphi}{\partial t}(\tilde{\boldsymbol\varphi}-\boldsymbol\varphi) dx \!\!\! 
&+&\!\!\!\gamma\varepsilon\int_{\Omega}\nabla\boldsymbol\varphi:\nabla(\tilde{\boldsymbol\varphi}-\boldsymbol\varphi)
dx+\frac{\gamma}{\varepsilon}\int_{\Omega}\Psi_0'(\boldsymbol\varphi)\cdot
(\tilde{\boldsymbol\varphi}
-\boldsymbol\varphi)dx\nonumber\\
&&\hspace{-3cm}-\frac12\beta J_{0}(\boldsymbol u,\boldsymbol\varphi)^{-1}\int_{\Omega}(\tilde{\varphi}^{N}-\varphi^{N})\, c \, |\boldsymbol u-\boldsymbol u_{\Omega}|^{2}\nonumber\\
&&\hspace{-3cm}-\int_{\Omega}(\tilde{\varphi}^{N}-\varphi^{N})\boldsymbol f\cdot(\alpha\boldsymbol u+\boldsymbol p)-\langle{\mathcal E}(\boldsymbol p),{\mathcal E}(\boldsymbol u)\rangle_{\mathbb C'(\boldsymbol\varphi)
(\tilde{\boldsymbol\varphi}-\boldsymbol\varphi)}\geq 0.\label{eq:gf}
\end{eqnarray}
In addition, $\boldsymbol u$ and $\boldsymbol p$ have to solve the 
state equation (SE) and the adjoint equation (AE), see Theorem
\ref{OPTSYS}.
The constraints $\varphi^N=0$ on $S_0$ and $\varphi^N=1$ on $S_1$ can be easily 
incorporated by imposing these conditions when a mesh point lies in 
$S_0\cup S_1$. 
We replace $\frac{\partial \boldsymbol\varphi}{\partial t}$ in (\ref{eq:gf}) by a time discrete approximation 
which corresponds to a pseudo time stepping approach. 
We then discretize the resulting inequality,   
the state equation (SE) and the adjoint equation (AE) 
using standard finite element approximations, see \cite{our}.

In the computation we present we take the weighting factor 
$c= 2000$ in $\Omega:=(-1,1)\times(-1,1)$ 
and $\mathbf{u}_\Omega=\mathbf{0}$. 
We set $\Gamma_D = \{(-1,y)\cup(1,y)\in \mathbb{R}^2: y \in [-1,-0.9]\cup[0.9,1]\}$ and 
$\Gamma_{g} = \Gamma_{g_-}\cup \Gamma_{g_+}$ 
with $\Gamma_{g_\pm} := \{(\pm 1,y)\in \mathbb{R}^2: y \in [-0.8,-0.7]\cup[-0.1,0.1]\cup[0.7,0.8]\}$. 
We take $\boldsymbol{g} =(\pm 7,0)^T$  on $\Gamma_{g_\pm}$ and $S_1=\emptyset$. 
Since we wish to have material 
adjacent to the parts of the boundary that are fixed and where the forces are applied we set 
$S_0=\{(x,y)\in \mathbb{R}^2: x \in [-1,-0.9]\cup[0.9,1], 
y \in [-1,-0.9]\cup[-0.8,-0.7]\cup[-0.1,0.1]\cup[0.7,0.8]\cup[0.9,1]\}$.  
We take $N=3$ and use an isotropic elasticity tensor
$\mathbb C^1$ of the form
$\mathbb C^1\mathcal{E}=2\mu_1\mathcal{E}+\lambda_1(tr\mathcal{E})I$ with
$\lambda_1=\mu_1 =10$
and we choose $\mathbb C^2 =\frac12 \mathbb C^1$ and 
$\mathbb C^3 = \varepsilon^2 \mathbb C^1$ in the void.  
The interfacial parameters we use are $\varepsilon = \tfrac1{18\pi}$ and 
$\gamma= 0.2$
and we set $\beta=10$. In addition, we choose the masses 
$\boldsymbol m = (0.35, 0.15, 0.5)^T$.

\begin{figure}[h]
\includegraphics[width = 4.5cm]{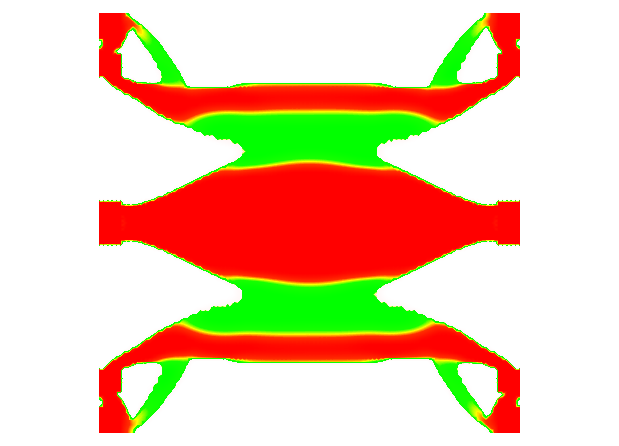}
\includegraphics[width = 4.5cm]{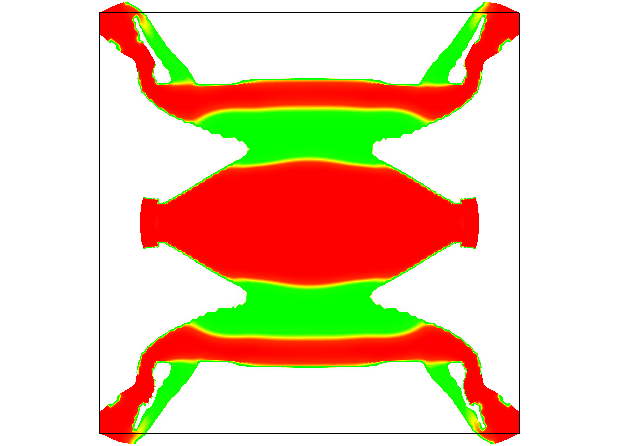}
\caption{Push simulation with three phases (left) and deformed configuration with the outline of the initial geometry (right)}
\label{push_rid}
\end{figure}

\begin{figure}[h]
\includegraphics[width = 4.5cm]{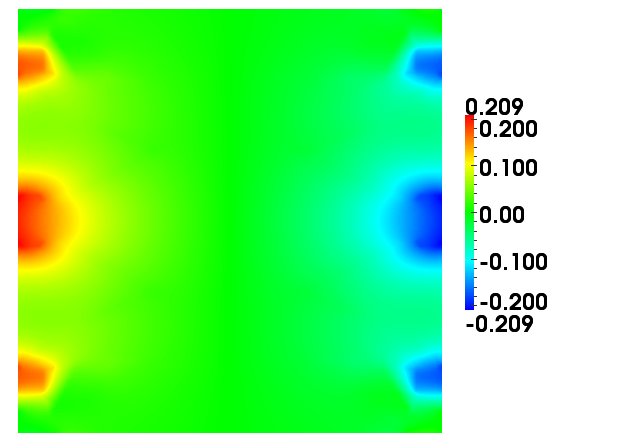} 
\includegraphics[width = 4.5cm]{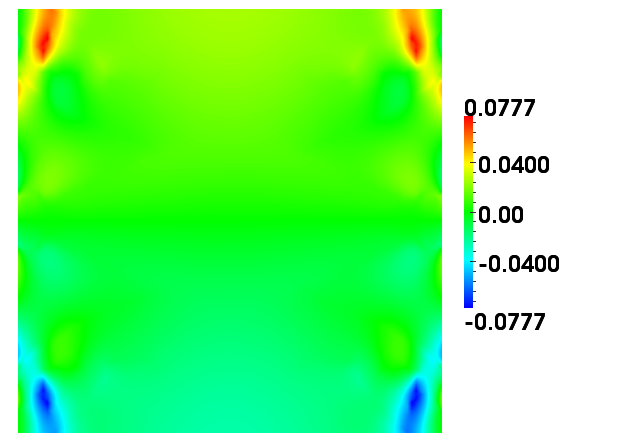} 
\caption{Displacement vector $\mathbf{u}$, $x$-component (left), $y$ component (right).}
\label{push_rid_dsp}
\end{figure}

In the Figure \ref{push_rid} we display the optimized configuration (left hand plot) 
and the deformed optimal configuration together with the outline of the initial geometry (right hand plot), here 
hard material is shown in red and soft material in green. 
In Figure \ref{push_rid_dsp}  we display the displacement vector $\mathbf{u}$. 


\end{document}